\documentclass[preprint,12pt]{elsarticle}

\usepackage [latin1]{inputenc}
\usepackage{amssymb}
\usepackage{amsfonts}
\usepackage{mathrsfs}
\usepackage{amsmath}
\usepackage{amscd}
\usepackage{amssymb}
\usepackage{amsmath}
\usepackage{graphicx}
\usepackage{array}
\usepackage{bm}

\journal{}
\def\proof{\par{\it Proof}. \ignorespaces}

\begin{document}

\newtheorem{prop}{Proposition}[section]
\newdefinition{remark}{Remark}[section]
\newtheorem{theorem}{Theorem}[section]
\newtheorem{lemma}{Lemma}[section]
\newtheorem{corollary}{Corollary}[section]
\newtheorem{defn}[theorem]{Definition}
\renewcommand{\theequation}{\thesection.\arabic{equation}}

\begin{frontmatter}



\title{Ergodicity for Stochastic Neutral Retarded Partial Differential Equations Driven
by $\alpha$-regular Volterra process}


\author{Xia Pan} \corref{cor1}\cortext[cor1]{Corresponding author.}\ead{panxia13@mails.ucas.ac.cn}

\address{College of Liberal Arts and Sciences, National University of Defense Technology, Changsha, 410073, P.R. China.}

\author{Zhi Li} \ead{}

\address{School of Mathematics and Information Sciences, Yangtze University, Jingzhou 434023, Hubei,
People¡¯s Republic of China}

\begin{abstract}
In this article, we study the ergodicity of neutral retarded stochastic functional differential equations driven by $\alpha$-regular Volterra process. Based on the equivalence between neutral retarded stochastic functional differential equations and the stochastic evolution equation, we get the ergodicity of neutral retarded stochastic functional differential equations.
\end{abstract}

\begin{keyword}
$\alpha$-regular Volterra process, neutral retarded stochastic functional differential equations, stochastic evolution equation, ergodicity




\end{keyword}

\end{frontmatter}



\section{Introduction}
In P. Coupek\cite{cou}, the authors considered the following stochastic evolution equation
\begin{equation}\label{1}
\left\{\begin{array}{l}
\mathrm{d} X_{t}=A X_{t}+\Phi \mathrm{d} B_{t}, \quad t \geq 0, \\
X_{0}=x,
\end{array}\right.
\end{equation}
where A generates a $C_{0}$-semigroup of bounded linear operators $S=(S(t), t \geq 0)$ acting on a separable Hilbert space and its mild solution is defined by
$$
X_{t}^{x}:=S(t) x+\int_{0}^{t} S(t-r) \Phi \mathrm{d} B_{r}, \quad t \geq 0.
$$
The driving process is a two-sided Hilbert space valued $\alpha$-regular Volterra process $B$ and $\Phi$ is a bounded linear operator. It is shown that the the solution of the stochastic evolution equation (\ref{1}) is a stationary process under some conditions. But the ergodicity of the stochastic evolution equation (\ref{1}) does not considered. Li\cite{li} discuss a class of neutral retarded stochastic functional differential equations driven by a fractional Brownian motion on Hilbert spaces, the ergodicity of the strictly stationary solution and non-stationary solution is studied. Liu\cite{liu1,liu2,liu3,liu4} has done a series of work on the stationary of the stochastic retarded evolution equations.

Motivated by their work, the purpose of this article is to study the ergodicity of neutral retarded stochastic functional differential equations driven by $\alpha$-regular Volterra process:
\begin{equation}\label{3.1}
\left\{\begin{array}{l}\mathrm{d}\left(x(t)-D x_{t}\right)=A\left(x(t)-D x_{t}\right) \mathrm{d} t+F x_{t} \mathrm{~d} t+B \mathrm{d} B(t), \quad t>0, \\ x(0)=\phi_{0}+D \phi_{1}, \quad x_{0}=\phi_{1}, \quad t \in[-r, 0], \quad \phi=\left(\phi_{0}, \phi_{1}\right) \in \mathcal{H}.\end{array}\right.
\end{equation}
In this article, we will prove the equation (\ref{3.1}) is equivalent to the following equation:
\begin{equation}\label{3.2}
\left\{\begin{array}{l}
\mathrm{d} X(t)=\mathcal{A} X(t) \mathrm{d} t+\mathcal{B} \mathrm{d} B(t), \quad t>0, \\
X(0)=\phi=\left(\phi_{0}, \phi_{1}\right) \in \mathcal{H}.
\end{array}\right.
\end{equation}
Subsequently, the ergodicity behavior of stationary solution and non-stationary solution for the equation (\ref{3.1}) is also investigated.

The rest of the paper is arranged as followings. In section 2, we first introduce some preliminaries on $\alpha$-regular Volterra process. Section 3 is devoted to the study of the equivalence of the equation (\ref{3.1}) and the equation (\ref{3.2}), besides, the ergodicity behavior of the equation (\ref{3.1}) is considered.

\section{Preliminaries}

In this section, we develop a $C_{0}$-semigroup theory\cite{cou,li,mao} of the driving deterministic neutral system and collect some notions, conceptions and lemmas on Wiener integrals with respect to two-sided $\alpha$-regular Volterra processes which will be used throughout the whole of this paper.

\subsection{Strongly Continuous Semigroups}

 Let $V$ be a separable Hilbert space and $a: V \times V \rightarrow \mathbb{R}$ be a bounded bilinear form satisfying
$$
a(x, x) \leq-\alpha\|x\|_{V}^{2}, \quad \forall x \in V,
$$
where $\alpha>0$. A linear operator $A$ is defined with the form
$$a(x, y)=\langle x, A y\rangle_{V, V^{*}}, \quad x, y \in V,$$
where $V^{*}$ is the dual space of $V$. Then, $A$ generates a $C_{0}$-semigroup $e^{t A},~ t \geq 0,$ on $V^{*}$.

Define
$$
H=\left\{x \in V^{*}: \int_{0}^{\infty}\left\|A e^{t A} x\right\|_{V^{*}}^{2} \mathrm{d} t<\infty\right\}=\left(V, V^{*}\right)_{1 / 2,2},
$$
with inner product
$$
\langle x, y\rangle_{H}=\langle x, y\rangle_{V^{*}}+\int_{0}^{\infty}\left\langle A e^{t A} x, A e^{t A} y\right\rangle_{V^{*}} \mathrm{d} t, \quad x, y \in V^{*}.
$$
Denote the dual of $H$ by $H^{*}$, then $V \hookrightarrow H=H^{*} \hookrightarrow V^{*}$ and $\|x\|_{H}^{2} \leq \beta\|x\|_{V}^{2},~x \in V,$ for some $\beta>0$, where the imbedding $\hookrightarrow$ is dense and continuous with for some constant $\beta>0$. Moreover, for any $T \geq 0$ it is well known that
$$
L^{2}([0, T] ; V) \cap W^{1,2}\left([0, T] ; V^{*}\right) \subset C([0, T] ; H),
$$
where $W^{1,2}\left([0, T] ; V^{*}\right)$ is the Sobolev space and $C([0, T] ; H)$ is the space of all continuous functions from $[0, T]$ into $H$. It can be also shown that the semigroup $e^{t A}, t \geq 0,$ is bounded and analytic on both $V^{*}$ and $H$ such that $e^{t A}: V^{*} \rightarrow V$ for each $t>0$ and for some constant $M_{0}>0,$
$$
\left\|e^{t A}\right\| \mathscr{L}\left(V^{*}\right) \leq M_{0}, \quad\left\|e^{t A}\right\|_{\mathscr{L}(H)} \leq e^{-\alpha t} \text { for all } t \geq 0.
$$

Throughout the paper,
$$L_{r}^{2}=L^{2}([-r, 0] ; V)=\{\varphi(\theta):\int_{-r}^{0}\|\varphi(\theta)\|_{V}^{2} \mathrm{~d} \theta<\infty,r>0\}.$$
$\mathcal{H}=H \times L_{r}^{2}$ with its respective norm and inner product defined by
$$
\|\Phi\|_{\mathcal{H}}=\left(\left\|\phi_{0}\right\|_{H}^{2}+\left\|\phi_{1}\right\|_{L_{r}^{2}}^{2}\right)^{\frac{1}{2}}, \quad\langle\Phi, \Psi\rangle_{\mathcal{H}}=\left\langle\phi_{0}, \psi_{0}\right\rangle_{H}+\left\langle\phi_{1}, \psi_{1}\right\rangle_{L_{r}^{2}},
$$
where $\Phi=\left(\phi_{0}, \phi_{1}\right), \Psi=\left(\psi_{0}, \psi_{1}\right) \in \mathcal{H} .$

Let $x_{t}(\theta):=x(t+\theta)$ for any $t \geq 0$ and $\theta \in[-r, 0]$.

Suppose that $D_{1} \in \mathscr{L}(V)$, $\left.\left.D_{2} \in \mathscr{L}(L^{2}([-r, 0]) ; V\right) ; V\right)$, $F_{1} \in \mathscr{L}\left(V, V^{*}\right)$, $F_{2} \in \mathscr{L}\left(L^{2}([-r, 0]) ; V\right) ; V^{*})$, and two linear mappings $D$ and $F$ on $C([-r, T] ; V),$ respectively, denoted  by
$$
D x_{t}=D_{1} x(t-r)+D_{2} x_{t}, \quad t \in[0, T], \quad \forall x(\cdot) \in C([-r, T] ; V),
$$
and
$$
F x_{t}=F_{1} x(t-r)+F_{2} x_{t}, \quad t \in[0, T], \quad \forall x(\cdot) \in C([-r, T] ; V),
$$

Now, we consider the following deterministic functional differential equation of neutral type in $V^{*}$,
\begin{equation}
\label{1.1}
\left\{\begin{array}{l}
\frac{\mathrm{d}}{\mathrm{d} t}\left(x(t)-D x_{t}\right)=A\left(x(t)-D x_{t}\right)+F x_{t}, \quad \text { for any } \quad t \geq 0, \\
x(0)=\phi_{0}+D \phi_{1}, \quad x_{0}=\phi_{1}, \quad t \in[-r, 0], \quad \phi=\left(\phi_{0}, \phi_{1}\right) \in \mathcal{H}.
\end{array}\right.
\end{equation}
The integral form of (\ref{1.1}) is given by
$$
\left\{\begin{array}{l}
x(t)-D x_{t}=e^{t A} \phi_{0}+\int_{0}^{t} e^{(t-s) A} F x_{s} \mathrm{~d} s, \quad \text { for any } t \geq 0, \\
x(0)=\phi_{0}+D \phi_{1}, \quad x_{0}=\phi_{1}, \quad t \in[-r, 0], \quad \phi=\left(\phi_{0}, \phi_{1}\right)\in\mathcal{H}.
\end{array}\right.
$$
We say that $x$ is a (strict) solution of (\ref{1.1}) in $[0, T]$ if $x \in L^{2}([0, T] ; V) \cap$ $W^{1,2}\left([0, T] ; V^{*}\right)$ and the equation (\ref{1.1}) is satisfied almost everywhere in $[0, T]$ $T \geq 0$
Let $x(t), t \geq-r$ denote the unique solution of system (\ref{1.1}) with initial $x(0)=\phi_{0}+$ $D \phi_{1}$ and $x_{0}=\phi_{1}, \phi=\left(\phi_{0}, \phi_{1}\right) \in \mathcal{H} .$ We define a family of operators $\mathcal{S}(t): \mathcal{H} \rightarrow \mathcal{H}$
$t \geq 0,$ by
$$
\mathcal{S}(t) \phi=\left(x(t)-D x_{t}, x_{t}\right), \text { for any } \phi \in \mathcal{H}.
$$
It can be shown that the mapping $\mathcal{S}(t), t \geq 0,$ is a $C_{0}$ -semigroup with its infinitesimal generator $\mathcal{A}$ on the space $\mathcal{H}$.

\begin{lemma}  The family $t \rightarrow \mathcal{S}(t)$ is a strongly continuous semigroup on $\mathcal{H}$, i.e.,
\begin{itemize}
\item $\mathcal{S}(t) \in \mathscr{L}(\mathcal{H})$ for each $t \geq 0$,
\item $\mathcal{S}(0)=I, \mathcal{S}(t+s)=\mathcal{S}(t) \mathcal{S}(s)$ for any $s, t \geq 0$,
\item $\lim _{t \rightarrow 0^{+}} \mathcal{S}(t) \phi=\phi$ for each $\phi \in \mathcal{H}$.
\end{itemize}
\end{lemma}

Moreover, the generator $\mathcal{A}$ may be explicitly specified as follows.

\begin{lemma}\label{2.2} The generator $\mathcal{A}$ of the strongly continuous semigroup $\mathcal{S}(t), t \geq 0,$ is described by
$$
\mathscr{D}(\mathcal{A})=\left\{\left(\phi_{0}, \phi_{1}\right) \in \mathcal{H}: \phi_{1} \in W^{1,2}([-r, 0] ; V), \phi_{0}=\phi_{1}(0)-D \phi_{1} \in V, A \phi_{0}+F \phi_{1} \in H\right\}
$$
and for each $\phi=\left(\phi_{0}, \phi_{1}\right) \in \mathscr{D}(\mathcal{A}),$
$$
\mathcal{A} \phi=\left(A \phi_{0}+F \phi_{1}, \phi_{1}^{\prime}\right) \in \mathcal{H}.
$$

\end{lemma}

\subsection{Two-sided Volterra processes}

Let $K: \mathbb{R}^{2} \rightarrow \mathbb{R}$ be a kernel such that: $K(t, r)=0$ on $\{t<r\}$ and $\lim _{t \rightarrow r+} K(t, r)=0$ for every $r \in \mathbb{R}$, and $K(\cdot, r)$ is continuously differentiable in $(r, \infty)$ for every $r \in \mathbb{R}$, besides, there is an $\alpha \in\left(0, \frac{1}{2}\right)$ such that
$$
\left|\frac{\partial K}{\partial u}(u, r)\right| \lesssim(u-r)^{\alpha-1}
$$
on $\{r:r<u\},$ where $A \lesssim B$ means that there is a finite positive constant $C$ such that $A \leq C B$, The constant $C$ is independent of all the changeable arguments of the expressions $A$ and $B$.
Such a function K is called an $\alpha$-regular Volterra kernel.

Let $K$ be an $\alpha$-regular Volterra kernel, then
$$\phi(u, v):=\int_{-\infty}^{u \wedge v} \frac{\partial K}{\partial u}(u, r) \frac{\partial K}{\partial v}(v, r) \mathrm{d} r \lesssim|u-v|^{2 \alpha-1}.
$$

Define
$$
R\left(s_{1}, t_{1}, s_{2}, t_{2}\right):=\int_{\mathbb{R}}\left(K\left(t_{1}, r\right)-K\left(s_{1}, r\right)\right)\left(K\left(t_{2}, r\right)-K\left(s_{2}, r\right)\right) \mathrm{d} r,
$$
then, for $s_{1}<t_{1}$ and $s_{2}<t_{2}$
\begin{equation}\label{2.1}
R\left(s_{1}, t_{1}, s_{2}, t_{2}\right)=\int_{s_{1}}^{t_{1}} \int_{s_{2}}^{t_{2}} \phi(u, v) \mathrm{d} u \mathrm{d} v.
\end{equation}

\begin{defn} A stochastic process $b=\left(b_{t}, t \in \mathbb{R}\right)$ is an $\alpha$-regular Volterra process if it is centered $b_{0}=0$ a.s. and such that
$$
\mathbb{E}\left(b_{t_{1}}-b_{s_{1}}\right)\left(b_{t_{2}}-b_{s_{2}}\right)=R\left(s_{1}, t_{1}, s_{2}, t_{2}\right),
$$
for every $s_{1}, s_{2}, t_{2}, t_{2} \in \mathbb{R},$ where $R$ is defined by formula (\ref{2.1}) with an $\alpha$ -regular Volterra kernel $K$.
\end{defn}

\subsection{Wiener integration}
Let $\left(V,\langle\cdot, \cdot\rangle_{V}\right)$ be a separable Hilbert space. Let $b=\left(b_{t}, t \in \mathbb{R}\right)$ be a two-sided Volterra process with a kernel $K$. Denote by $\mathscr{E}(\mathbb{R} ; V)$ the set of $V$ -valued step functions on $\mathbb{R}$, i.e. $f \in \mathscr{E}(\mathbb{R} ; V)$ satisfies
$$
f=\sum_{j=1}^{n} f_{j} \mathbf{1}_{\left[t_{j-1}, t_{j}\right)},
$$
where $n\in\mathbb{N}, -\infty<t_{0}<t_{1}<\cdots<t_{n}<\infty$ and $f_{j} \in V$ for all $j=1,2, \ldots, n .$ Consider the linear mapping $i: \mathscr{E}(\mathbb{R} ; V) \rightarrow$ $L^{2}(\Omega ; V)$ given by
$$
i: \quad f:=\sum_{j=1}^{n} f_{i} \mathbf{1}_{\left[t_{j-1}, t_{j}\right)} \longmapsto \sum_{j} f_{j}\left(b_{t_{j}}-b_{t_{j-1}}\right)=: i(f)
$$
and define the operator $\mathscr{K}^{*}: \mathscr{E}(\mathbb{R} ; V) \rightarrow L^{2}(\mathbb{R} ; V)$ by
$$
\left(\mathscr{K}^{*} f\right)(r):=\int_{r}^{\infty} f(u) \frac{\partial K}{\partial u}(u, r) \mathrm{d} u, \quad r \in \mathbb{R}.
$$
For simplicity, it is assumed here that $\mathscr{K}^{*}$ is injective.
We have
$$
\|i(f)\|_{L^{2}(\Omega ; V)}=\left\|\mathscr{K}^{*} f\right\|_{L^{2}(\mathbb{R} ; V)}
$$
for $f \in \mathscr{E}(\mathbb{R} ; V)$. Now, have $\mathscr{E}(\mathbb{R} ; V)$ completed under the inner product
$$
\langle f, g\rangle_{\mathscr{D}}:=\left\langle\mathscr{K}^{*} f, \mathscr{K}^{*} g\right\rangle_{L^{2}(\mathbb{R}, V)}
$$
denote the completion by $\mathscr{D}(\mathbb{R} ; V)$ and extend $\mathscr{K}^{*}$ to $\left(\mathscr{D}(\mathbb{R} ; V),\langle\cdot, \cdot\rangle_{\mathscr{D}}\right)$ which is now a
Hilbert space. This in turn extends $i$ to a linear isometry between $\mathscr{D}(\mathbb{R} ; V)$ and a closed linear subspace of $L^{2}(\Omega ; V)$. The space $\mathscr{D}(\mathbb{R} ; V)$ is viewed as the space of admissible integrands and, for $f \in \mathscr{D}(\mathbb{R} ; V),$ the random variable $i(f)$ is the stochastic integral of $f$ with respect to the Volterra process $b$.

\begin{defn} Let $U$ be a real separable Hilbert space. An $\alpha$ -regular $U$ -cylindrical Volterra process is a collection $B=\left(B_{t}, t \in \mathbb{R}\right)$ of bounded linear operators $B_{t}: U \rightarrow L^{2}(\Omega)$ such that
- for every $u \in U, B(u)$ is a centered stochastic process in $\mathbb{R}$ with $B_{0}(u)=0$ a.s.;
- for every $s_{1}, t_{1}, s_{2}, t_{2} \in \mathbb{R}$ and every $u_{1}, u_{2} \in U$ it holds that
$$
\mathbb{E}\left(B_{t_{1}}\left(u_{1}\right)-B_{s_{1}}\left(u_{2}\right)\right)\left(B_{t_{2}}\left(u_{2}\right)-B_{s_{2}}\left(u_{2}\right)\right)=R\left(s_{1}, t_{1}, s_{2}, t_{2}\right)\left\langle u_{1}, u_{2}\right\rangle_{V}
$$
with $R$ given by formula (\ref{2.1}).
\end{defn}

\section{Ergodic Theorems}
Let $U, V$ be two real separable Hilbert spaces and consider the stochastic evolution equation. In this section, we will investigate the ergodic theorem for the solution to
\begin{equation}\label{31}
\left\{\begin{array}{l}\mathrm{d}\left(x(t)-D x_{t}\right)=A\left(x(t)-D x_{t}\right) \mathrm{d} t+F x_{t} \mathrm{~d} t+B \mathrm{d} B(t), \quad t>0, \\ x(0)=\phi_{0}+D \phi_{1}, \quad x_{0}=\phi_{1}, \quad t \in[-r, 0], \quad \phi=\left(\phi_{0}, \phi_{1}\right) \in \mathcal{H}\end{array}\right.
\end{equation}
where $A$ is an infinitesimal generator of a strongly continuous semigroup $(S(t), t \geq 0)$ of bounded linear operators acting on $V$ and $x \in L^{2}(\Omega ; V)$. We assume that $\Phi \in \mathscr{L}(U, V)$ and $B=\left(B_{t}, t \in \mathbb{R}\right)$ is an $\alpha$ -regular $U$ -cylindrical Volterra process. We can associate (3.1) with an abstract stochastic differential equation without delay and neutral item on $\mathcal{H}$,
\begin{equation}\label{32}
\left\{\begin{array}{l}
\mathrm{d} X(t)=\mathcal{A} X(t) \mathrm{d} t+\mathcal{B} \mathrm{d} B(t), \quad t>0,\\
X(0)=\phi=\left(\phi_{0}, \phi_{1}\right) \in \mathcal{H}
\end{array}\right.
\end{equation}
where $\mathcal{A}$ is given as in Lemma \ref{2.2}, $\mathcal{B}: U \rightarrow \mathcal{H}$ is defined by
$$
\mathcal{B}: x \rightarrow(B x, 0), \quad \forall x \in U.
$$
The solution to equation (\ref{32}) is given in the mild form by the variation of constants formula
$$
X_{t}^{x}:=S(t) x+Z_{t}:=S(t) x+\int_{0}^{t} S(t-r) \mathcal{B} \mathrm{d} B_{r}, \quad t \geq 0.
$$

We shall show below that (\ref{31}) and (\ref{32}) are actually equivalent in the sense that every solution $t \rightarrow X(t)$ of (\ref{32}) is of the form
$$
X(t)=\left(x(t)-D x_{t}, x_{t}\right), \quad t \geq 0,
$$
where the function $x(t)$ is the solution of (\ref{31}) .

\begin{prop} For any $T \geq 0,$ let $x:[-r, T] \times \Omega \rightarrow H$ be a solution of (\ref{31}) Then the mapping $X:(t, \omega) \rightarrow\left(x(t, \omega)-D x_{t}(\omega), x_{t}(\omega)\right) \in \mathcal{H}$ from $[0, T] \times \Omega$ into
$\mathcal{H}$ is a solution of the abstract Cauchy problem (\ref{32}).
\end{prop}

\proof Recall that
$$
e^{t \mathcal{A}} \phi=\left(\tilde{x}(t)-D \tilde{x}_{t}, \tilde{x}_{t}\right), \quad t \in[0, T],
$$
where $\tilde{x}(t)$ is the unique solution of the homogeneous equation
\begin{equation}\label{3.3}
\left\{\begin{array}{l}
\mathrm{d}\left(\tilde{x}(t)-D \tilde{x}_{t}\right)=A\left(\tilde{x}(t)-D \tilde{x}_{t}\right) \mathrm{d} t+F \tilde{x}_{t} \mathrm{~d} t, \quad \text { for any } \quad t>0, \\
\tilde{x}(0)=\phi_{0}+D \phi_{1}, \quad \tilde{x}_{0}=\phi_{1}, \quad t \in[-r, 0], \quad \phi=\left(\phi_{0}, \phi_{1}\right) \in \mathcal{H}.
\end{array}\right.
\end{equation}
On the other hand, for any $h \in H,$ let $\phi_{0}=h, \phi_{1}=0$ for $\theta \in[-r, 0]$ and $\phi=(h, 0),$ we define the fundamental solution $G(t)$ of (\ref{3.3}) with this initial datum by(see\cite{li},\cite{huang})
$$
G(t) h=\left\{\begin{array}{ll}
x(t, \phi), & t \geq 0, \\
0, & t<0.
\end{array}\right.
$$
Then, by virtue of fundamental solution $G(t), t \in \mathbb{R}^{1},$ we have
$$
e^{t \mathcal{A}}\left(\phi_{0}, 0\right)=\left(G(t) \phi_{0}-D G_{t} \phi_{0}, G_{t} \phi_{0}\right), \quad \phi_{0} \in H,
$$
and
$$
x(t)=\tilde{x}(t)+\int_{0}^{t} G(t-s) B \mathrm{~d} B(s), \quad t \in[0, T],
$$
which further yields the equality
$$
x_{t}=\tilde{x}_{t}+\int_{0}^{t} G_{t-s} B \mathrm{~d} B(s), \quad t \in[0, T].
$$
Indeed, for $\theta \in[-r, 0],$ the integral $\int_{0}^{t} G(t+\theta-s) B \mathrm{~d} B(s)$ is equal to $\int_{0}^{t+\theta} G(t-$
s) $B \mathrm{~d} B(s)$ if $t+\theta \geq 0,$ and to 0 if $t+\theta<0$ since $G(t+\theta-s)=0$ for $s \in[t+\theta, t] .$ Then, for any $t \in[0, T]$
$$
X(t)=e^{t \mathcal{A}} \phi+\int_{0}^{t} e^{(t-s) \mathcal{A}} \mathcal{B} \mathrm{d} B(s).
$$
That is, $X(t)$ is the solution of the abstract Cauchy problem (\ref{32}).
\qed

 \begin{prop}
 Let $X(t)=\left(\pi_{0} X(t), \pi_{1} X(t)\right) \in \mathcal{H}$ be the solution of (\ref{32}) with arbitrary initial $\phi=\left(\phi_{0}, \phi_{1}\right) \in \mathcal{H} .$ Then the process
$$
x(t)=\left\{\begin{array}{ll}
\pi_{0} X(t)+D\left(\pi_{0} X\right)_{t}, & \text { if } t \geq 0, \\
x(0)=\phi_{0}+D \phi_{1}, & x_{0}=\phi_{1}, \text { if } t \in[-r, 0]\\
\end{array}\right.
$$
is a solution of (\ref{31}) and $\pi_{1} X(t)=\left(\pi_{0} X\right)_{t}$ for all $t \geq 0$.
\end{prop}
\proof Since $X$ is a solution of the problem
$$
\left(\mathrm{d} \pi_{0} X(t), \mathrm{d} \pi_{1} X(t)\right)=\left(A \pi_{0} X(t)+F\left(\pi_{1} X\right)(t), \mathrm{d} \pi_{1} X(t)\right) \mathrm{d} t+\left(B\mathrm{d} B(t), 0\right),
$$
it, thus, follows that for all $t \geq 0$,
$$
\pi_{0} X(t)=e^{t A} \phi_{0}+\int_{0}^{t} e^{(t-s) A} F\left(\pi_{1} X\right)(s) \mathrm{d} s+\int_{0}^{t} e^{(t-s) A} B \mathrm{~d} B(s).
$$

On the other hand, one can notice by definition that
$$
\left(\pi_{0} X\right)_{t}=\left\{\begin{array}{ll}
\left(\pi_{0} X\right)(t+\theta), & \text { for } t+\theta \geq 0,\\
\phi_{1}(t+\theta), & \text { for } t+\theta<0,
\end{array}\right.
$$
where $\phi_{1} \in L^{2}([-r, 0] \times \Omega ; H)$.

Therefore,
$$
\left(\pi_{1} X\right)(t)=\left(\pi_{0} X\right)_{t} \in L^{2}([-r, 0] \times \Omega ; H)
$$
for all $t \geq 0$ and $x(t)$ is a mild solution of (\ref{31}).
\qed

According to P.Coupek\cite{cou}, the equation (\ref{32}) has stationary solution.

(H) Let $S(r) \Phi \in \mathscr{L}_{2}(U, V)$ for all $r>0$. Let there further exist $T_{0}>0$ such that
$$
\int_{0}^{T_{0}}\|S(r) \Phi\|_{\mathscr{L}_{2}(U, V)}^{\frac{2}{1+2 \alpha}} \mathrm{d} r<\infty.
$$

\begin{defn} Let $d \geq 1$. We say that an $\mathbb{R}^{d}$ -valued stochastic process $Y=\left(Y_{t}, t \in \mathbb{R}\right)$ has
- stationary increments if for every $n \in \mathbb{N}$ and every $s_{i}, t_{i} \in \mathbb{R}, s_{i}<t_{i}, i=1,2, \ldots, n,$ we have that the following holds for every $h \in \mathbb{R}$ :
$$
\begin{array}{l}
\operatorname{Law}\left(Y_{t_{1}}-Y_{s_{1}}, Y_{t_{2}}-Y_{s_{2}}, \ldots, Y_{t_{n}}-Y_{s_{n}}\right) \\
\quad=\operatorname{Law}\left(Y_{t_{1}+h}-Y_{s_{1}+h}, Y_{t_{2}+h}-Y_{s_{2}+h}, \ldots, Y_{t_{n}+h}-Y_{s_{n}+h}\right),
\end{array}
$$
reflexive increments if for every $n \in \mathbb{N}$ and every $s_{i}, t_{i} \in \mathbb{R}, s_{i}<t_{i}, i=1,2, \ldots, n,$ we have that the following holds:
$$
\begin{array}{l}
\operatorname{Law}\left(Y_{t_{1}}-Y_{s_{1}}, Y_{t_{2}}-Y_{s_{2}}, \ldots, Y_{t_{n}}-Y_{s_{n}}\right) \\
\quad=\operatorname{Law}\left(Y_{-s_{1}}-Y_{-t_{1}}, Y_{-s_{2}}-Y_{-t_{2}}, \ldots, Y_{-s_{n}}-Y_{-t_{n}}\right).
\end{array}
$$
\end{defn}

\begin{lemma}\cite{cou} Assume that $B$ has stationary and reflexive increments. Assume further that $S(u) \Phi \in \mathscr{L}_{2}(U, V)$ for every $u>0$ and that
$$
\int_{0}^{\infty}\|S(r) \Phi\|_{\mathscr{L}_{2}(U, V)}^{\frac{2}{1+2 \alpha}} \mathrm{d} r<\infty
$$
holds. Then there is a measure $\mu_{\infty}$ such that
$$
w^{*}-\lim _{t \rightarrow \infty} \mu_{t}^{0}=\mu_{\infty}.
$$
\end{lemma}

\begin{theorem}
(Ergodic theorem for a stationary solution)

Let $\left(x^{\tilde{x}_{0}+D \tilde{x}_{1}}(t), t \geq 0\right)$ be a $\mathcal{H}$-value stationary solution to (\ref{31}). Let $\varrho: H \rightarrow \mathbb{R}$ be a measurable functional such that $\mathbb{E}\left|\varrho\left(\tilde{x}_{0}+D \tilde{x}_{1}\right)\right|<\infty$. Then
$$\lim _{T \rightarrow \infty} \frac{1}{T} \int_{0}^{T} \varrho\left(x^{\tilde{x}_{0}+D \tilde{x}_{1}}(t)\right) \mathrm{d} t=\int_{H} \varrho(y) \mu_{\infty}^{*}(\mathrm{~d} y), \quad a.s. -\mathbb{P}.$$
\end{theorem}

\proof we know that there exists $\tilde{x}=\left(\tilde{x}_{0}, \tilde{x}_{1}\right) \in \mathcal{H}$, a random variable on $(\Omega, \mathscr{F}, \mathbb{P})$ such that $X^{\tilde{x}}(t)=\left(X_{1}(t), X_{2}(t)\right)$ is a
stationary solution to (\ref{32}) and $X^{\tilde{x}}(t)=\left(X_{1}(t), X_{2}(t)\right)$ is ergodic. On the other hand,
$$x(t)=\left\{\begin{array}{ll}X_{1}(t)+D\left(X_{1}\right)_{t}, & \text { if } ~t \geq 0, \\ x(0)=\tilde{x}_{0}+D \tilde{x}_{1}, & x_{0}=\tilde{x}_{1}, \quad \text { if }~ t \in[-r, 0]\end{array}\right.$$
is the solution of (\ref{31}). Notice that $\left(X_{1}(t), X_{2}(t)\right)$ is a stationary solution to (\ref{32})
and $\left(X_{1}(t), X_{2}(t)\right)$ is ergodic, then $X_{1}(t)$ and $\left(X_{1}\right)_{t}$ are two stationary processes and
$X_{1}(t)$ and $\left(X_{1}\right)_{t}$ are ergodic. Thus, we have $x(t)=X_{1}(t)+D\left(X_{1}\right)_{t}, t>0$ with $x(0)=\tilde{x}_{0}+D \tilde{x}_{1}$ also is ergodic.
\qed

Now, we are in a position to consider the ergodic theorem for an arbitrary solution
to (\ref{31}).

\begin{theorem}
If the semigroup $\mathcal{S}(t)$ is exponentially stable, i.e., there exist constants
$M>0$ and $\rho>0$ such that for all $t \geq 0,\|\mathcal{S}(t)\| \mathscr{L}(\mathcal{H}) \leq M e^{-\rho l}$. Let $\left(x^{x_{0}}(t), t \geq 0\right)$ be a solution to (\ref{31}) and $\varrho: H \rightarrow \mathbb{R}$ be a functional satisfying the global Lipschitz condition, i.e., there exists a constant $L > 0$ such that
$$|\varrho(x)-\varrho(y)| \leq L\|x-y\|_{H},$$
for all $x, y \in H$. Then
$$\lim _{T \rightarrow \infty} \frac{1}{T} \int_{0}^{T} \varrho\left(x^{x_{0}}(t)\right) \mathrm{d} t=\int_{H} \varrho(y) \mu_{\infty}^{*}(\mathrm{~d} y), \quad a.s. -\mathbb{P},$$
for all $x_{0} \in H$.
\end{theorem}

\proof The desired convergence can be rewritten as
$$\lim _{T \rightarrow \infty}\left|\frac{1}{T} \int_{0}^{T} \varrho\left(x^{x_{0}}(t)\right) \mathrm{d} t-\int_{H} \varrho(y) \mu_{\infty}^{*}(\mathrm{~d} y)\right|=0, \quad a.s. -\mathbb{P},$$
for all $x_{0} \in H .$

Let $\left(x^{\tilde{x}_{0}+D \tilde{x}_{1}}(t), t \geq 0\right)$ be a $H$ -value stationary solution to $(\ref{31}) .$ Then
$$\begin{aligned}
& \left|\frac{1}{T}\int_{0}^{T} \varrho\left(x^{x_{0}}(t)\right) \mathrm{d} t-\int_{H} \varrho(y) \mu_{\infty}^{*}(\mathrm{~d} y)\right| \\
 \leq &\mid \frac{1}{T} \int_{0}^{T} \varrho\left(x^{x_{0}}(t)\right) \mathrm{d} t-\frac{1}{T} \int_{0}^{T} \varrho\left(x^{\tilde{x}_{0}+D \tilde{x}_{1}}(t)\right) \\
 &+ \left|\frac{1}{T} \int_{0}^{T} \varrho\left(x^{\tilde{x}_{0}+D \tilde{x}_{1}}(t)\right) \mathrm{d} t-\int_{H} \varrho(y) \mu_{\infty}^{*}(\mathrm{~d} y)\right|\\
=:& I_{1}+I_{2}\end{aligned}.$$

Using the Lipschitz assumption and the exponential stability of $S(t),$ we get
$$\begin{aligned} I_{1} &=\left|\frac{1}{T} \int_{0}^{T} \varrho\left(x^{x_{0}}(t)\right) \mathrm{d} t-\frac{1}{T} \int_{0}^{T} \varrho\left(x^{\tilde{x}_{0}+D \tilde{x}_{1}}(t)\right) \mathrm{d} t\right| \\ & \leq \frac{L}{T} \int_{0}^{T}\left\|x^{x_{0}}(t)-x^{\tilde{x}_{0}+D \tilde{x}_{1}}(t)\right\|_{H} \mathrm{~d} t \\ &=\frac{L}{T} \int_{0}^{T}\left\|S(t)\left(x_{0}-\tilde{x}_{0}-D \tilde{x}_{1}\right)\right\|_{H} \mathrm{~d} t \\ & \leq \frac{L}{T}\left\|x_{0}-\tilde{x}_{0}-D \tilde{x}_{1}\right\|_{H} \int_{0}^{T} e^{-\alpha t} \mathrm{~d} t, \end{aligned}$$
which goes to zero as $T \rightarrow \infty$. On the other hand, by Ergodic theorem for a stationary solution,  $I_{2}$ goes to zero as $T \rightarrow \infty$.
\qed

\noindent{\bf Acknowledgements}: This paper is supported by  Key Laboratory of Engineering mathematical modeling and analysis of Hunan Province (CN)(2019MMAED002).

\end{document}